\documentclass[11pt]{amsart}
\usepackage{amsmath}
\usepackage{amssymb}
\usepackage{amscd}
\usepackage{url}
\usepackage[latin2]{inputenc}
\usepackage[T1]{fontenc}
\usepackage{graphicx}

\def\Ent{{\mathrm {Ent}}}

\def\E{{\mathbb E}}
\def\R{{\mathbb R}}
\def\PP{{\mathbb P}}

\def\1{{\mathbf 1}}
\def\ud{{\mathrm d}}
\def\e{{\mathrm e}}

\theoremstyle{definition}
\newtheorem{example}{Example}[section]

\theoremstyle{plain}

\newtheorem{lemma}[example]{Lemma}

\newtheorem{thm}[example]{Theorem}

\newtheorem{corollary}[example]{Corollary}

\theoremstyle{remark}

\newtheorem{remark}[example]{Remark}
\numberwithin{equation}{section}

\title[Boolean functions with small second order influences]{Boolean functions with small second order influences
on the discrete cube} 
\author{Krzysztof Oleszkiewicz} 
\address{Institute of Mathematics, University of Warsaw, ul. Banacha 2,
02-097 Warsaw, Poland}
\email{koles@mimuw.edu.pl}
\thanks{Research supported by the National Science Centre, Poland, project number 2015/18/A/ST1/00553.}

\date{}

\begin{document}
\begin{abstract}
Motivated by a recent paper of Kevin Tanguy~\cite{KT}, in which the concept of second order influences on the discrete cube and Gauss space has been investigated in detail, the present note studies it in a more specific context of Boolean functions on the discrete cube. Some bounds which Tanguy obtained in~\cite{KT} as applications of his more general approach are extended and complemented. 
\end{abstract}
\maketitle

\section{Introduction}

Throughout the paper, $n$ stands for an integer greater than $1$, and we use the standard notation $[n]:=\{ 1, 2, \ldots, n\}$. We 
equip the discrete cube $\{ -1, 1\}^{n} = \{-1,1\}^{[n]}$ with the normalized counting (
uniform probability) measure
$\mu_{n}=({\frac{1}{2}\delta_{-1}+\frac{1}{2}\delta_{1}})^{\otimes n}$. Let $\E$
denote the expectation with respect to this measure,
and let $r_{1},$ $r_{2}, \ldots,$ $r_{n}$ be the standard Rademacher functions on the discrete cube, i.e. the coordinate projections $r_{i}(x)=x_{i}$
for $x \in \{ -1,1\}^{n}$ and $i \in [n]$. Furthermore, for $A \subseteq [n]$, we define the Walsh functions by $w_{A}=\prod_{i \in A}r_{i}$, with
$w_{\emptyset} \equiv 1$.

The Walsh functions $(w_{A})_{A \subseteq [n]}$ form a complete orthonormal system in $L^{2}\left(\{-1,1\}^{n},\mu_{n}\right)$. Thus, every
$f:\{-1,1\}^{n} \to \R$ admits a unique Walsh-Fourier expansion
$
f=\sum_{A \subseteq [n]} {\hat{f}}(A)w_{A},
$
whose coefficients are given by
\[
{\hat{f}}(A)=\langle f, w_{A} \rangle = \E[f \cdot w_{A}]=\frac{1}{2^{n}}\,\sum_{x \in \{-1,1\}^{n}} f(x)w_{A}(x).
\]

For $p \geq 1$, we abbreviate $\|f\|_{L^{p}(\{-1,1\}^{n},\mu_{n})}$ to $\| f\|_{p}$.

In a standard way, for $i \in [n]$ we define the influence of the $i$-th coordinate on the function $f$
by
\[
I_{i}=I_{i}(f):=\sum_{A \subseteq [n]:\,i \in A} \big({\hat{f}}(A)\big)^{2}=\sum_{A \subseteq [n]:\,i \in A} {\hat{f}}^{2}(A),
\]
and given distinct integers $i,j \in [n]$, we define (following, up to a minor modification, the notation of Kevin Tanguy's paper~\cite{KT}) the influence of the couple $(i,j)$ on the function~$f$ by
\[
I_{i,j}=I_{i,j}(f):=\sum_{A \subseteq [n]:\,i,j \in A} {\hat{f}}^{2}(A).
\]
Note that if $f$ is Boolean, i.e.~$\{-1,1\}$-valued, then $I_{i}, I_{i,j} \leq 1$ since, by the orthonormality of the Walsh-Fourier system,
$\sum_{A \subseteq [n]} {\hat{f}}^{2}(A)=\E\big[f^{2}\big]=1$.

\section{Main results} \label{druga}

\begin{thm} \label{main} There exists a universal constant $C>0$ with the following property. Let $n \geq 2$ be an integer. Assume that 
$f:\,\{-1,1\}^{n} \to \{-1,1\}$
satisfies the bound
$I_{i,j}(f) \leq \alpha n^{-2}\ln^{2}n$ for all $1 \leq i<j \leq n$ for some $\alpha>0$. Then
${\hat{f}}^{2}(\emptyset) \geq 1-C\alpha$, or there exists exactly one $i \in [n]$ such that
${\hat{f}}^{2}(\{i\}) \geq 1-C\alpha n^{-1}\ln n$.
\end{thm}

We postpone the proof of Theorem~\ref{main} till 
Section~\ref{dowod}. This theorem says that if a Boolean function $f$ on the discrete cube has uniformly small second order influences, then it has to be close to one of the functions $1, -1$, or very close to one of the dictatorship/antidictatorship functions $r_{1}, -r_{1}, \ldots, r_{n}, -r_{n}$. Thus, it may be viewed as a modified version of two classical theorems: the KKL theorem of Kahn, Kalai, and Linial,~\cite{KKL}, which says that if a Boolean function $f$ on the discrete cube has influences $I_{i}$ uniformly bounded from above by $\alpha n^{-1}\ln n$, then ${\hat{f}}^{2}(\emptyset) \geq 1-C\alpha$, where $C>0$ is some universal constant, and the FKN theorem of Friedgut, Kalai, and Naor,~\cite{FKN}, which says that if a Boolean function on the discrete cube is close -- in a certain sense, different from the assumptions of Theorem~\ref{main} -- to an affine function, then it must be close to one of the constant functions or to one of the dictatorship/antidictatorship functions.

While the numerical value of the constant $C$ which can be deduced from the proof of Theorem~\ref{main} is quite large (with some additional effort, many numerical constants in the proofs can be improved, though perhaps at the cost of clarity), it is not difficult to obtain a reasonable estimate in the case of $\alpha$ close to zero.

\begin{thm} \label{main2} There exists a bounded function $C:\,(0,\infty) \to (0,\infty)$ such that 
$\lim_{\alpha \to 0^{+}} C(\alpha)=4$ and with
the following property. Let $n \geq 2$ be an integer. Assume that a function
$f:\,\{-1,1\}^{n} \to \{-1,1\}$
satisfies the bound
$\max_{1 \leq i <j \leq n} I_{i,j}(f) \leq \alpha n^{-2}\ln^{2}n$ 
for some $\alpha>0$. Then
${\hat{f}}^{2}(\emptyset) \geq 1-C(\alpha)\alpha$, or there exists exactly one $i \in [n]$ such that
${\hat{f}}^{2}(\{i\}) \geq 1-C(\alpha)\alpha n^{-1}\ln n$.
\end{thm}

\section{Auxiliary notation and tools}

For any $g:\,\{-1,1\}^{n} \to \R$ and $t\geq 0$, we define $P_{t}g:\,\{-1,1\}^{n} \to \R$ by
\[
P_{t}g=\sum_{A \subseteq [n]} \e^{-|A|t}\,{\hat{g}}(A)w_{A}.
\]
$(P_{t})_{t \geq 0}$ is called the heat semigroup on the discrete cube. By the classical hypercontractive inequality of Bonami~\cite{B}, $\| P_{t}g\|_{2} \leq \| g\|_{1+\e^{-2t}}$.

Given $f:\,\{-1,1\}^{n} \to \R$ and $i \in [n]$, define the discrete partial derivative $D_{i}f:\,\{-1,1\}^{n} \to \R$ of~$f$ by $(D_{i}f)(x)=\big(f(x)-f(x^{i})\big)/2$, where
\[
x^{i}=(x_{1},\ldots, x_{i-1},-x_{i},x_{i+1},\ldots,x_{n})\,\,\,\hbox{for}\,\,\,x \in \{-1,1\}^{n}. 
\]
Given distinct integers $i,j \in [n]$, let $D_{i,j}=D_{i} \circ D_{j}$, so that $D_{i,j}f=D_{i}(D_{j}f)$.
Note that
\[
D_{i}f=\sum_{A \subseteq [n]:\,i \in A} {\hat{f}}(A)w_{A}\,\,\,\hbox{and}\,\,\,D_{i,j}f=\sum_{A \subseteq [n]:\,i,j \in A} {\hat{f}}(A)w_{A},
\]
so that $I_{i}(f)=\| D_{i}f\|_{2}^{2}$ and $I_{i,j}(f)=\| D_{i,j}f\|_{2}^{2}$.

A slightly different partial derivative operator $\partial_{i}$ is defined by the formula
$(\partial_{i}f)(x)=\big(f(x^{i \to 1})-f(x^{i \to \,-1})\big)/2$, where 
\[
x^{i \to \varepsilon}=(x_{1},\ldots, x_{i-1},\varepsilon,x_{i+1},\ldots,x_{n})\,\,\,\hbox{for}\,\,\,x \in \{-1,1\}^{n}, \,\,\varepsilon \in \{-1,1\}.
\]
Again, given distinct integers $i,j \in [n]$, let $\partial_{i,j}=\partial_{i} \circ \partial_{j}$, i.e. $\partial_{i,j}f=\partial_{i}(\partial_{j}f)$.
One easily checks that $D_{i}f=r_{i} \cdot \partial_{i}f$ and thus $D_{i,j}f=r_{i}r_{j} \cdot \partial_{i,j}f$ for all functions~$f$. The Rademacher functions $r_{i}$ and $r_{j}$ are $\{-1,1\}$-valued, so $\| D_{i}f\|_{p}=\| \partial_{i}f\|_{p}$ and
$\| D_{i,j}f\|_{p}=\| \partial_{i,j}f\|_{p}$ for every $p \geq 1$. Furthermore, $P_{t}D_{i,j}f=\e^{-2t}r_{i}r_{j} \cdot P_{t}\partial_{i,j}f$ (this identity is an easy consequence of the same equality for Walsh functions which is easy to verify) and thus $\e^{2t}\,\| P_{t}D_{i,j}f\|_{p}=\|P_{t}\partial_{i,j}f\|_{p}$ for all $p \geq 1$, $t \geq 0$, and $f:\,\{-1,1\}^{n} \to \R$.

\begin{lemma} \label{above1}
For every $f:\,\{-1,1\}^{n} \to \R$ there is
\[
\sum_{A \subseteq [n]:\,|A| \geq 2} {\hat{f}}^{2}(A)=
4\,\sum_{i,j:\,1 \leq i<j\leq n} \int_{0}^{\infty} \big(\e^{2t}-1\big)\|P_{t}D_{i,j}f\|_{2}^{2}\,\ud t.
\]
\end{lemma}

\begin{proof} 
Since $P_{t}D_{i,j}f=\sum_{A \subseteq [n]:\,i,j \in A} \e^{-|A|t}{\hat{f}}(A)w_{A}$, we have
\[
\sum_{i,j:\, i<j} \|P_{t}D_{i,j}f\|_{2}^{2} = \sum_{i,j:\,i<j} \,\sum_{A \subseteq [n]:\,i,j \in A} \e^{-2|A|t}
\,{\hat{f}}^{2}(A)
\]
\[
=\sum_{A \subseteq [n]:\,|A| \geq 2} \,\, \sum_{i,j \in A:\,i<j} \e^{-2|A|t}\,{\hat{f}}^{2}(A)
=\sum_{A \subseteq [n]:\,|A| \geq 2} {\,|A|\, \choose 2}  \e^{-2|A|t}\,{\hat{f}}^{2}(A), 
\]
so that
\[
\sum_{i,j:\,1 \leq i<j\leq n} \int_{0}^{\infty} \big(\e^{2t}-1\big)\|P_{t}D_{i,j}f\|_{2}^{2}\,\ud t
\]
\[
=\sum_{A \subseteq [n]:\,|A| \geq 2} \frac{|A|(|A|-1)}{2}\,{\hat{f}}^{2}(A) \cdot 
\int_{0}^{\infty} \big(\e^{2t}-1\big)\e^{-2|A|t}\,\ud t.
\]
It remains to note that, for $k>1$,
\[
\int_{0}^{\infty} \big(\e^{2t}-1\big)\e^{-2kt}\,\ud t=\frac{1}{2(k-1)}-\,\frac{1}{2k}=\frac{1}{2k(k-1)}.
\]
\end{proof}

\begin{lemma} \label{estimate}
For every $f:\,\{-1,1\}^{n} \to \{-1,1\}$ and integers $1 \leq i<j \leq n$ there is
\[
\int_{0}^{\infty} \big(\e^{2t}-1\big)\|P_{t}D_{i,j}f\|_{2}^{2}\,\ud t \leq \frac{2I_{i,j}}{\ln^{2}(2/I_{i,j})},
\]
where $I_{i,j}$ denotes the influence of the pair $(i,j)$ on the function $f$.
\end{lemma}

\begin{proof}
Bonami's hypercontractive bound applied to $g=\partial_{i,j}f$ yields
\[
\e^{4t}\| P_{t}D_{i,j}f\|_{2}^{2}=
\| P_{t}\partial_{i,j}f\|_{2}^{2} \leq
\| \partial_{i,j}f\|_{1+\e^{-2t}}^{2}
=\big(\E\big[|\partial_{i,j}f|^{1+\e^{-2t}}\big]\big)^{\frac{2}{1+\e^{-2t}}}.
\]
Note that $\partial_{i,j}f$ is $\{-1, -1/2,0,1/2,1\}$-valued because $f$ is Boolean. Since
$|w|^{1+\e^{-2t}} \leq 2^{1-\e^{-2t}}\cdot w^{2}$ for every $w \in \{-1, -1/2,0,1/2,1\}$, we also have
\[
\E\big[|\partial_{i,j}f|^{1+\e^{-2t}}\big] \leq 2^{1-\e^{-2t}}\cdot \E\big[(\partial_{i,j}f)^{2}\big]=2^{1-\e^{-2t}}\cdot \|\partial_{i,j}f\|_{2}^{2}=2^{1-\e^{-2t}} \cdot I_{i,j},
\]
therefore
\[
\| P_{t}D_{i,j}f\|_{2}^{2} \leq \e^{-4t} \cdot 4^{\frac{1-\e^{-2t}}{1+\e^{-2t}}} \cdot I_{i,j}^{\frac{2}{1+\e^{-2t}}}
=\frac{(1-u)^{2}}{(1+u)^{2}} \cdot 4^{u} \cdot I_{i,j}^{1+u},
\]
where we use the change of variables $u=\frac{1-\e^{-2t}}{1+\e^{-2t}} \in [0,1]$, i.e. $t=\frac{\ln(1+u)-\ln(1-u)}{2}$,
and thus $\frac{\ud t}{\ud u}=\frac{1}{2}\left(\frac{1}{1+u}+\frac{1}{1-u}\right)=\frac{1}{(1-u)(1+u)}$. Hence
\[
\int_{0}^{\infty} \big(\e^{2t}-1\big)\|P_{t}D_{i,j}f\|_{2}^{2}\,\ud t \leq \int_{0}^{1} \frac{2u}{1-u} \cdot \frac{(1-u)^{2}}{(1+u)^{2}} \cdot 4^{u} \cdot I_{i,j}^{1+u} \cdot \frac{\ud u}{(1-u)(1+u)}
\]
\[
=2\int_{0}^{1} \frac{4^{u}\,u}{(1+u)^{3}}I_{i,j}^{1+u}\,\ud u  \leq 
2I_{i,j}\int_{0}^{\infty} u\left(\frac{I_{i,j}}{2}\right)^{u}\,\ud u=\frac{2I_{i,j}}{\ln^{2}(2/I_{i,j})}.
\]
We have used the fact that, by the convexity of the exponential function, $2^{u} \leq 1+u$ for $u \in [0,1]$.
\end{proof}

\begin{lemma} \label{m_i}
For every $f:\,\{-1,1\}^{n} \to \R$ and $i \in [n]$ there is
\[
I_{i}(f)-{\hat{f}}^{2}(\{i\})=
\sum_{A \subseteq [n]:\,|A| \geq 2,\, i \in A} {\hat{f}}^{2}(A)=
2\,\sum_{j \in [n]\setminus \{i\}} \int_{0}^{\infty} \e^{2t}\,\|P_{t}D_{i,j}f\|_{2}^{2}\,\ud t.
\]
\end{lemma}

\begin{proof} 
Since $P_{t}D_{i,j}f=\sum_{A \subseteq [n]:\,i,j \in A} \e^{-|A|t}{\hat{f}}(A)w_{A}$, we have
\[
\sum_{j \in [n] \setminus \{i\}} \|P_{t}D_{i,j}f\|_{2}^{2}
= \sum_{j \in [n] \setminus \{i\}} \,\sum_{A \subseteq [n]:\,i,j \in A} \e^{-2|A|t}\,{\hat{f}}^{2}(A)
\]
\[
=\!\sum_{A \subseteq [n]:\,|A| \geq 2,\,i \in A} \,\, \sum_{j \in A \setminus \{i\}} \e^{-2|A|t}\,{\hat{f}}^{2}(A)
=\!\sum_{A \subseteq [n]:\,|A| \geq 2,\,i \in A} \!\!(|A|-1)\,\e^{-2|A|t}\,{\hat{f}}^{2}(A), 
\]
so that
\[
\sum_{j \in [n]\setminus \{i\}} \int_{0}^{\infty} \e^{2t}\,\|P_{t}D_{i,j}f\|_{2}^{2}\,\ud t
=\sum_{A \subseteq [n]:\,|A| \geq 2,\,i \in A} \!\!(|A|-1)\,{\hat{f}}^{2}(A)
\int_{0}^{\infty} \!\e^{2t}\,\e^{-2|A|t}\,\ud t
\]
\[
=\sum_{A \subseteq [n]:\,|A| \geq 2,\,i \in A} (|A|-1)\,{\hat{f}}^{2}(A)\,\frac{1}{2|A|-2}=
\frac{1}{2}\sum_{A \subseteq [n]:\,|A| \geq 2,\,i \in A} {\hat{f}}^{2}(A).
\]
\end{proof}

\begin{lemma} \label{estimate2}
For every $f:\,\{-1,1\}^{n} \to \{-1,1\}$ and integers $1 \leq i<j \leq n$ there is
\[
\int_{0}^{\infty} \e^{2t}\,\|P_{t}D_{i,j}f\|_{2}^{2}\,\ud t \leq \frac{I_{i,j}}{\ln(1/I_{i,j})},
\]
where $I_{i,j}$ denotes the influence of the pair $(i,j)$ on the function $f$.
\end{lemma}

\begin{proof}
In the same way as in the proof of Lemma~\ref{estimate}, we obtain the bound
\[
\| P_{t}D_{i,j}f\|_{2}^{2} \leq \e^{-4t} \cdot 4^{\frac{1-\e^{-2t}}{1+\e^{-2t}}} \cdot I_{i,j}^{\frac{2}{1+\e^{-2t}}}
=\frac{(1-u)^{2}}{(1+u)^{2}} \cdot 4^{u} \cdot I_{i,j}^{1+u},
\]
where again $u=\frac{1-\e^{-2t}}{1+\e^{-2t}} \in [0,1]$, $t=\frac{\ln(1+u)-\ln(1-u)}{2}$,
and $\frac{\ud t}{\ud u}=\frac{1}{(1-u)(1+u)}$. Hence
\[
\int_{0}^{\infty} \e^{2t}\,\|P_{t}D_{i,j}f\|_{2}^{2}\,\ud t \leq \int_{0}^{1} \frac{1+u}{1-u} \cdot \frac{(1-u)^{2}}{(1+u)^{2}} \cdot 4^{u} \cdot I_{i,j}^{1+u} \cdot \frac{\ud u}{(1-u)(1+u)}
\]
\[
=\int_{0}^{1} \frac{4^{u}}{(1+u)^{2}}\,I_{i,j}^{1+u}\,\ud u  \leq 
\int_{0}^{\infty} I_{i,j}^{1+u}\,\ud u = \frac{I_{i,j}}{\ln(1/I_{i,j})}.
\]
We have once more used the fact that $2^{u} \leq 1+u$ for $u \in [0,1]$.
\end{proof}

\begin{lemma} \label{elementary}
Let $z \in [0,1/4)$ and $0 \leq x \leq y \leq x^{2}+z$. Then either $y \leq 2z$ or $x \geq 1-2z$.
\end{lemma}

\begin{proof}
Note that $\big(\frac{1}{4}-z\big)^{1/2} \geq \frac{1}{2}-2z$ for $z \in [0,1/4]$, so it suffices to prove
that $y \leq x_{1}$ or $x \geq x_{2}$, where
\[
x_{1}= \frac{1}{2}-\Big(\frac{1}{4}-z\Big)^{1/2}\,\,\,\mathrm{and}\,\,\,x_{2}= \frac{1}{2}+\Big(\frac{1}{4}-z\Big)^{1/2}.
\]
This is easy: since $x \leq x^{2}+z$, we have $x \in (-\infty, x_{1}] \cup [x_{2},\infty)$. Thus, if $x < x_{2}$, then $x \in [0, x_{1}]$, which in turn implies that
\[
y \leq x^{2}+z \leq x_{1}^{2}+z=x_{1}.
\]
\end{proof}

We will also make use of the following slightly more precise observation.

\begin{remark} \label{asymptota}
Let $z \in [0,1/4)$ and $0 \leq x \leq y \leq x^{2}+z$. Then either $y \leq z+4z^{2}$ or $x \geq 1-z-4z^{2}$. To prove this, it suffices to replace in the proof of Lemma~\ref{elementary} the $\big(\frac{1}{4}-z\big)^{1/2} \geq \frac{1}{2}-2z$ bound by a stronger one, $\big(\frac{1}{4}-z\big)^{1/2} \geq \frac{1}{2}-z-4z^{2}$, which
is also satisfied for all $z \in [0,1/4]$.
\end{remark}

\begin{lemma} \label{trick}
For every $f:\,\{-1,1\}^{n} \to \{-1,1\}$ and $i \in [n]$, the influence $I_{i}$ of the $i$-th coordinate on the function $f$ satsfies the inequality
\[
|{\hat{f}}(\{i\})| \leq I_{i}=|{\hat{f}}(\{i\})|^{2}+\sum_{A \subseteq [n]:\,|A| \geq 2,\, i \in A} {\hat{f}}^{2}(A).
\]
\end{lemma}

\begin{proof}
Using the triangle inequality and the fact that $\partial_{i}f$ is $\{-1,0,1\}$-valued,
so that~$|\partial_{i}f| \equiv (\partial_{i}f)^{2}$, we arrive at
\[
|{\hat{f}}(\{i\})|=|\E[\partial_{i}f]| \leq \E[|\partial_{i}f|]=\E\big[(\partial_{i}f)^{2}\big]=I_{i}.
\]
\end{proof}

\section{Proof of the main results} \label{dowod}

\noindent{\em Proof of Theorem~\ref{main}.} 
Let us define a positive constant $\kappa$ by
\[
\kappa=\min\left(\inf_{n \geq 2} \frac{n}{16 \ln n},\,\, \inf_{n \geq 2} \frac{n}{\ln^{2}n}\right),
\]
and let $C=\max\big(\kappa^{-1}, 20\big)$.
For $\alpha \geq \kappa$ the assertion of the theorem holds true in a trivial way, therefore we may and will
assume that $0<\alpha<\kappa$.

Since $\alpha < n/\ln^{2}n$ and $I_{i,j} \leq \alpha n^{-2}\ln^{2}n$, we have 
\[
\ln(2/I_{i,j}) \geq \ln(1/I_{i,j}) > \ln n,
\]
and thus also
\[
\frac{I_{i,j}}{\ln^{2}(2/I_{i,j})} \leq \frac{\alpha}{n^{2}}
\,\,\,\,\,\mathrm{and}\,\,\,\,\,
\frac{I_{i,j}}{\ln(1/I_{i,j})} \leq \frac{\alpha \ln n}{n^{2}}
\]
for all $1 \leq i<j \leq n$, so that
\[
\sum_{i,j:\,1 \leq i <j \leq n} \frac{I_{i,j}}{\ln^{2}(2/I_{i,j})} \leq {n \choose 2} \cdot \frac{\alpha}{n^{2}} \leq \frac{\alpha}{2}
\]
and, for every $i \in [n]$,
\[
\sum_{j \in [n] \setminus \{i\}} \frac{I_{i,j}}{\ln(1/I_{i,j})} \leq (n-1) \cdot \frac{\alpha \ln n}{n^{2}} \leq \frac{\alpha \ln n}{n}.
\]
Using Lemma~\ref{above1} and Lemma~\ref{estimate}, we arrive at
\[
\sum_{A \subseteq [n]:\,|A| \geq 2} {\hat{f}}^{2}(A) \leq 4\alpha,
\]
and from Lemma~\ref{m_i} and Lemma~\ref{estimate2} we obtain, for every $i \in [n]$,
\[
\sum_{A \subseteq [n]:\,|A| \geq 2,\,i \in A} {\hat{f}}^{2}(A) \leq 2\alpha n^{-1}\ln n. 
\]
Applying Lemma~\ref{elementary} to $x=|{\hat{f}}(\{i\})|$, $y=I_{i}(f)$, and
$z=2\alpha n^{-1}\ln n<1/4$, we deduce from Lemma~\ref{trick} that, for every $i \in [n]$, either $I_{i}(f) \leq 4\alpha n^{-1}\ln n$ or $|{\hat{f}}(\{i\})| \geq 1-4\alpha n^{-1}\ln n$.

\medskip
\noindent{\bf Case 1: $I_{i}(f) \leq 4\alpha n^{-1}\ln n$ for all $i \in [n]$.} Then, by Lemma~\ref{trick},
for all $i \in [n]$ we have
$
|{\hat{f}}(\{i\})| \leq 4\alpha n^{-1}\ln n
$,
so that
\[
\sum_{i \in [n]} {\hat{f}}^{2}(\{i\}) \leq n \cdot 16\alpha^{2} n^{-2}\ln^{2}n < 16\alpha
\]
and therefore
\[
{\hat{f}}^{2}(\emptyset)=1-\sum_{A \subseteq [n]:\,|A|=1} {\hat{f}}^{2}(A)-\sum_{A \subseteq [n]:\,|A|\geq 2} {\hat{f}}^{2}(A) \geq 1-16\alpha-4\alpha \geq 1-C\alpha.
\]

\medskip
\noindent{\bf Case 2: $|{\hat{f}}(\{i\})| \geq 1-4\alpha n^{-1}\ln n$ for some $i \in [n]$.} Then
we obviously have ${\hat{f}}^{2}(\{i\}) \geq 1-8\alpha n^{-1}\ln n>1/2$, which implies that there is exactly one such
$i \in [n]$ (recall that $\sum_{A \subseteq [n]} {\hat{f}}^{2}(A)$=1). It remains
to note that $C>8. \,\,\,\Box$

\medskip
\noindent{\em Proof of Theorem~\ref{main2}.} It suffices to repeat the proof of Theorem~\ref{main}, using
Remark~\ref{asymptota} instead of Lemma~\ref{elementary} when $\alpha$ is close to zero.

\rightline{$\Box$}

\begin{remark}
Throughout the paper, we have restricted our interest to the uniform estimate assumption. However, it is easy to see that the assumption is used in the proof only via the
\[
\alpha \geq \max\left(\sum_{i,j:\,1 \leq i <j \leq n} \frac{2I_{i,j}}{\ln^{2}(2/I_{i,j})},\,\, \frac{n}{\ln n}\cdot \max_{i \in [n]} \sum_{j \in [n] \setminus \{i\}} \frac{I_{i,j}}{\ln(1/I_{i,j})}\right)
\]
condition, allowing for significant extensions. In particular, if a function $f:\,\{-1,1\}^{n} \to \{-1,1\}$ satisfies
\[
\max_{1 \leq i<j \leq n} I_{i,j}(f) \leq n^{-\gamma}\,\,\,\mathrm{and}\,\,\,\max_{i \in [n]} \sum_{j \in [n] \setminus \{i\}} I_{i,j}(f) \leq \beta \,\frac{\ln^{2}n}{n}
\]
for some constants $\beta>0$ and $\gamma \in (0,1]$, then the proof works for
$\alpha = \beta\gamma^{-2}$,\\
so that ${\hat{f}}^{2}(\emptyset) \geq 1-C\beta\gamma^{-2}$, or there exists a unique $i \in [n]$ such that
${\hat{f}}^{2}(\{i\}) \geq 1-C\beta\gamma^{-2}n^{-1}\ln n$, where $C>0$ is the universal constant from Theorem~\ref{main}.
\end{remark}

\section{Alternative proof}

Here we present another proof of Theorem~\ref{main}, based on a quite different approach. We will prove the following more general statement.

\begin{thm} \label{additional}
Let $n \geq 2$.
For a function $f:\,\{-1,1\}^{n} \to \{-1,1\}$, let
\[
\theta=\max_{i \in [n]} \,\sum_{j \in [n] \setminus \{i\}} I_{i,j}(f).
\]
Assume $\theta \leq 1/25$. Then either there exists exactly one $i \in [n]$ such that
\[
\big|\hat{f}(\{i\})\big| \geq 1-\frac{8\theta}{\ln(1/\theta)}
\]
or
\[
\max_{i \in [n]} \,I_{i}(f) \leq \frac{4\theta}{\ln(1/\theta)}.
\]
\end{thm}

For a function $\varphi:\,\{-1,1\}^{n} \to \R$, let us denote by $I(\varphi)$ its total influence,
\[
I(\varphi)=\sum_{i \in [n]} I_{i}(\varphi).
\]
Since $P_{0}\varphi \equiv \varphi$, Bonami's bound $\forall_{t \geq 0}\,\,\,\| P_{t}\varphi\|_{2} \leq 
\| \varphi\|_{1+\e^{-2t}}$ yields
\[
\frac{\ud}{\ud t} \| P_{t}\varphi\|_{2}\Big|_{t=0^{+}} \leq \frac{\ud}{\ud t}
\| \varphi\|_{1+\e^{-2t}}\Big|_{t=0^{+}},
\]
which amounts to the classical logarithmic Sobolev inequality,
\[
2\,I(\varphi) \geq \Ent\big(\varphi^{2}\big)=\E\big[\varphi^{2}\ln\big(\varphi^{2}\big)\big]-
\E\big[\varphi^{2}\big]\,\ln\E\big[\varphi^{2}\big].
\]
If the function $\varphi$ is $\{0,1\}$-valued, then $\varphi^{2} \equiv \varphi$
and $\varphi^{2}\ln\big(\varphi^{2}\big) \equiv 0$. In this well-known way, we derive a weak functional version of the edge-isoperimetric inequality:

\begin{lemma} \label{isoperimetry}
For every function $h:\,\{-1,1\}^{n} \to \{0,1\}$ there is
\[
I(h) \geq \frac{1}{2} \cdot \E[h] \cdot \ln(1/\E[h])=\frac{1}{2} \cdot \PP(h=1) \cdot \ln\big(1/\PP(h=1)\big).
\]
\end{lemma}

\begin{remark} A slightly stronger version of this inequality is known to hold true, with $\ln$ replaced by $\log_{2}$. It can be deduced from Harper's solution \cite{H} of the edge-isoperimetric problem for the discrete cube (this result has also been proved by Lindsey, Bernstein, and Hart).
\end{remark}

\begin{corollary} \label{3values}
For each $n \geq 2$
and every function $g:\,\{-1,1\}^{n} \to \{-1,0,1\}$ satisfying $I(g) \leq 1/25$ there exists $\eta \in \{-1,0,1\}$
such that
\[
\PP(g \neq \eta) \leq 4 \cdot I(g)/\ln\big(1/I(g)\big).
\]
\end{corollary}

\begin{proof} Since $\PP(g=-1)+\PP(g=0)+\PP(g=1)=1$, certainly there exists $\eta \in \{-1,0,1\}$
for which $\PP(g=\eta) \geq 1/3$, so that $\PP(g \neq \eta) \leq 2/3$. We will prove that the bound of Corollary~\ref{3values} holds true for this $\eta$. Let us define $h:\,\{-1,1\}^{n} \to \{0,1\}$
by $h=1_{g \neq \eta}$. Then, for every $x \in \{-1,1\}^{n}$ and every $i \in [n]$,
\[
|\partial_{i}h(x)|=\big|h\big(x^{i \to 1}\big)-h\big(x^{i \to -1}\big)\big|/2=
\big|1_{g \neq \eta}\big(x^{i \to 1}\big)-1_{g \neq \eta}\big(x^{i \to -1}\big)\big|/2
\]
\[
=
\frac{1}{2}\Big|1_{\{-1,0,1\}\setminus \{\eta\}}\Big(g\big(x^{i \to 1}\big)\Big)-1_{\{-1,0,1\}\setminus \{\eta\}}\Big(g\big(x^{i \to -1}\big)\Big)\Big|
\]
\[
\leq \big|g\big(x^{i \to 1}\big)-g\big(x^{i \to -1}\big)\big|/2=|\partial_{i}g(x)|,
\]
because $1_{\{-1,0,1\}\setminus\{\eta\}}:\,\{-1,0,1\} \to \{0,1\}$ is $1$-Lipschitz. Therefore
\[
I(h)=\sum_{i \in [n]} I_{i}(h)=\sum_{i \in [n]} \| \partial_{i}h\|_{2}^{2}
\leq \sum_{i \in [n]} \| \partial_{i}g\|_{2}^{2}=\sum_{i \in [n]} I_{i}(g)=I(g),
\]
so that, by Lemma~\ref{isoperimetry},
\[
\PP(h=1) \leq 2I(h)/\ln\big(1/\PP(h=1)\big) \leq 2I(g)/\ln\big(1/\PP(h=1)\big)
\]
\[
=2I(g)/\ln\big(1/\PP(g \neq \eta)\big) \leq
2I(g)/\ln(3/2) \leq 5I(g) \leq \sqrt{I(g)}.
\]
Thus, applying Lemma~\ref{isoperimetry} again, we get
\[
\PP(g \neq \eta)=\PP(h=1) \leq \frac{2I(h)}{\ln\big(1/\PP(h=1)\big)}
\leq \frac{2I(g)}{\ln\big(1/\sqrt{I(g)}\,\big)}=4 \cdot \frac{I(g)}{\ln\big(1/I(g)\big)}.
\]
\end{proof}

\noindent{\em Proof of Theorem~\ref{additional}.} Let $i \in [n]$. Since the function $f$ is
$\{-1,1\}$-valued, its $i$-th partial derivative $\partial_{i}f$ is $\{-1,0,1\}$-valued. Note that
\[
I(\partial_{i}f)=\sum_{j \in [n]} I_{j}(\partial_{i}f)=\sum_{j \in [n] \setminus \{i\}} \| \partial_{j}(\partial_{i}f)\|_{2}^{2}=\sum_{j \in [n] \setminus \{i\}} I_{i,j}(f) \leq \theta \leq 1/25.
\]
Applying Corollary~\ref{3values} to $g=\partial_{i}f$ and using the fact that the function $(0,1) \ni x \mapsto x/\ln(1/x)$ is increasing, we prove the existence of $\eta_{i} \in \{-1,0,1\}$
such that $\PP(\partial_{i}f \neq \eta_{i}) \leq 4\theta/\ln(1/\theta)$. Obviously, this $\eta_{i}$ is unique, because $4\theta/\ln(1/\theta)\leq 0.16/\ln(25)<1/20$, so that $\PP(\partial_{i}f=\eta_{i})>1/2$.

Observe that if $\partial_{i}f(x)=1$ for some $x=(x_{1},x_{2},\ldots,x_{n}) \in \{-1,1\}^{n}$, then $f\big(x^{i \to 1}\big)=1$ and $f\big(x^{i \to -1}\big)=-1$, i.e., $f(x)=x_{i}$. In other words, $\{ f \neq r_{i}\} \subseteq \{ \partial_{i}f \neq 1\}$. In a similar way, we prove 
$\{ f \neq -r_{i}\} \subseteq \{ \partial_{i}f \neq -1\}$. Thus, if there exists $i \in [n]$
for which $\eta_{i} \in \{-1,1\}$, then
\[
\PP(f \neq \eta_{i}r_{i}) \leq 4\theta/\ln(1/\theta),
\]
and since
\[
\hat{f}(\{i\})=\E[f \cdot r_{i}]=\PP(f=r_{i})-\PP(f \neq r_{i})=\PP(f \neq - r_{i})-\PP(f \neq r_{i})
\]
\[
=1-2\,\PP(f \neq r_{i})=2\,\PP(f \neq -r_{i})-1,
\]
we arrive at
\[
\big|\hat{f}(\{i\})\big| \geq 1-8\theta/\ln(1/\theta),
\]
as desired. The uniqueness of $i$ follows from $\hat{f}^{2}(\{i\}) > (9/10)^{2}>1/2$.

It remains to consider the case $\eta_{1}=\eta_{2}=\ldots=\eta_{n}=0$, in which
\[
I_{i}(f)=\E\big[(\partial_{i}f)^{2}\big]=\E[1_{\partial_{i}f \neq 0}]=\PP(\partial_{i}f \neq 0)
\leq 4\theta/\ln(1/\theta)
\]
for each $i \in [n]$, as desired. Since $\big|\hat{f}(\{i\})\big|\leq \sqrt{I_{i}(f)}$
and $4\theta/\ln(1/\theta)<1/20$, the assertion of Theorem~\ref{additional} is an exclusive disjunction, as stated.
$\,\,\,\Box$

\medskip

Theorem~\ref{main} easily follows from Theorem~\ref{additional} by considering $\theta=\alpha n^{-1}\ln^{2}n$ (and using the KKL theorem to derive the $\hat{f}^{2}(\emptyset) \geq 1-C\alpha$
estimate from the upper bound on $\max_{i \in [n]} I_{i}(f)$, as explained in Section~\ref{druga}). If $\alpha \leq 1/25$, then also $\theta \leq 1/25$, and by taking $C \geq 25$ we make the case $\alpha>1/25$ trivial.

\medskip

For the standard tribes function $T:\,\{-1,1\}^{n} \to \{-1,1\}$, with disjoint tribes of size $\sim \log_{2}(n/\ln n)$ each, one easily checks that $I_{i,j}(T) \simeq n^{-2}\ln^{2}n$ if $i$ and $j$ belong to different tribes, but $I_{i,j}(T) \simeq n^{-1}\ln n$ if $i$ and $j$ belong to the same tribe --~Remark (1) on page 703 of~\cite{KT} is not correct. Thus, $T$ cannot serve as an example showing the essential optimality of Theorem~\ref{main}, though to some extent it does the job for Theorem~\ref{additional}, with $\theta \simeq n^{-1}\ln^{2} n$ and $I_{i}(T) \simeq |{\hat{T}}(\{i\})|\simeq n^{-1}\ln n \simeq \theta/\ln(1/\theta)$ uniformly for all $i \in [n]$. 

This is complemented by the example of $V:\,\{-1,1\}^{n} \to \{-1,1\}$ defined by 
\[
V(x)=x_{1} \cdot \left(1-2\prod_{k=2}^{n}\Big(\frac{1+x_{k}}{2}\Big)\right).
\]
The function $V$ satisfies the assumptions of Theorem~\ref{additional} with $\theta \simeq n \cdot 2^{-n}$, i.e., $\theta/\ln(1/\theta) \simeq 2^{-n}$, whereas ${\hat{V}}(\{1\})=1-4 \cdot 2^{-n}$ and $I_{1}(V)=1$.

Very recently, Tomasz Przyby{\l}owski \cite{TP} established a nice counterpart of Theorem~\ref{main} for influences of orders higher than $2$. Following an advice of Peter Keevash, he also provided an example demonstrating the essential optimality of Theorem~\ref{main}.

\end{document}